\documentclass{amsart}
\usepackage{amscd}
\usepackage{amssymb}
\usepackage[ps,dvips,curve]{xypic}
\newbox\mybox
\def\overtag#1#2#3{\setbox\mybox\hbox{$#1$}\hbox to
  0pt{\vbox to 0pt{\vglue-#3\vglue-\ht\mybox\hbox to \wd\mybox
      {\hss$\ss#2$\hss}\vss}\hss}\box\mybox}
\def\undertag#1#2#3{\setbox\mybox\hbox{$#1$}\hbox to 0pt{\vbox to
    0pt{\vglue#3\vglue\ht\mybox\hbox to \wd\mybox
      {\hss$\ss#2$\hss}\vss}\hss}\box\mybox}
\def\lefttag#1#2#3{\hbox to 0pt{\vbox to 0pt{\vss\hbox to
      0pt{\hss$\ss#2$\hskip#3}\vss}}#1}
\def\righttag#1#2#3{\hbox to 0pt{\vbox to 0pt{\vss\hbox to
      0pt{\hskip#3$\ss#2$\hss}\vss}}#1}
\let\ss\scriptstyle

\def\Dot{\lower.2pc\hbox to 2pt{\hss$\bullet$\hss}}
\def\Circ{\lower.2pc\hbox to 2pt{\hss$\circ$\hss}}
\def\Vdots{\raise5pt\hbox{$\vdots$}}
\newcommand\lineto{\ar@{-}}
\newcommand\dashto{\ar@{--}}
\newcommand\dotto{\ar@{.}}

\newcommand\Ker{\operatorname{Ker}}
\newcommand\Hom{\operatorname{Hom}}
\newcommand\GL{\operatorname{GL}}
\newcommand\tr{\operatorname{tr}}
\newcommand\SL{\operatorname{SL}}

\newcommand\spec{\operatorname{Spec}}

\newcommand\Q{{\mathbb Q}}
\newcommand\C{{\mathbb C}}
\newcommand\Z{{\mathbb Z}}
\renewcommand\O{{\mathcal O}}
\newtheorem{theorem}{Theorem}[section]
\newtheorem*{theorem*}{Theorem}
\newtheorem{lemma}[theorem]{Lemma}
\newtheorem{proposition}[theorem]{Proposition}

\newtheorem{conjecture}[theorem]{Conjecture}
\newtheorem*{conjecture*}{Conjecture}

\theoremstyle{definition}

\newtheorem*{remark*}{Remark}

\begin{document}
\title{Universal abelian covers of quotient-cusps}
\author{Walter D. Neumann}
\address{Department of Mathematics\\Barnard College, Columbia
  University\\New York, NY 10027}
\email{neumann@math.columbia.edu}
\author{Jonathan Wahl}
\address{Department of Mathematics\\The University of
North Carolina\\Chapel Hill, NC 27599-3250}
\email{jw@math.unc.edu}
\thanks{This research was supported by grants from the Australian
  Research Council and the NSF (first author) and the 
  the NSA (second author).}
\begin{abstract}
  The quotient-cusp singularities are isolated complex surface
  singularities that are double-covered by cusp singularities.  We
  show that the universal abelian cover of such a singularity,
  branched only at the singular point, is a complete intersection cusp
  singularity of embedding dimension $4$. This supports a general
  conjecture that we make about the universal abelian cover of
  a $\Q$-Gorenstein singularity.
\end{abstract}
\maketitle

\section{Introduction}
\label{sec:intro}        

In this paper we introduce a general conjecture, and prove it in a new
case. Recall that a normal surface singularity is called
\emph{$\Q$-Gorenstein} if the dualizing sheaf $\omega_X$ has finite
order. Equivalently, for some $n>0$, $\Gamma(X-\{o\},
\omega_X^{\otimes n})$ is a free $\O_X$-module; or, $(X,o)$ is the
quotient of a Gorenstein singularity by a finite cyclic group acting
freely off the singular point.

\begin{conjecture}\label{conj:main}
  Let $(X,o)$ be a $\Q$-Gorenstein normal surface singularity whose
  link is a rational homology sphere.  Then the universal abelian
  cover $(\tilde X,o)$ of $(X,o)$ is a complete intersection, of
  embedding dimension at most the number of ends in the resolution
  graph of $(X,o)$.
\end{conjecture}

Our notation $(X,o)$ implies that we are considering the germ of a
normal complex surface singularity. Any finite connected covering of
the link of this singularity gives rise naturally to a unique analytic
covering of $(X,o)$ branched only at the singular
point (namely, form the cover $Z\to X-\{o\}$ and then take
$\spec\Gamma(\mathcal O_Z)$).  If the link is a rational homology sphere, then its
universal abelian covering is a finite covering, so it gives rise to
the so-called \emph{universal abelian covering} of the corresponding
singularity.

We shall discuss special cases of Conjecture \ref{conj:main} that
suggest that the complete intersection in the conjecture may be of
quite special type.  But first we discuss some of the background of
this conjecture.

For many years, the second author has been mystified by the
appropriate philosophical context by which to understand a theorem of
the other author, namely
\begin{theorem}[\cite{neumann83b}]\label{th:weighted homog}
Let $(X,o)$ be a weighted homogeneous normal surface singularity whose
link is a rational homology sphere.  Then the conclusion of Conjecture
\ref{conj:main} applies.  In fact, the universal abelian cover
$(\tilde X,o)$ of $(X,o)$ is a Brieskorn complete intersection.
\end{theorem}

The proof of the above theorem is completely convincing as well as
totally unenlightening: given the data determining the analytic type
of $(X,o)$, one writes down explicit equations of a Brieskorn complete
intersection plus an explicit action of an abelian group $G$, and then
argues directly that $G$ acts freely off the origin $o$, with quotient
isomorphic to $(X,o)$, and that $G$ has the correct order.

It is easy to see that $(X,o)$ is $\Q$-Gorenstein in the above
theorem, and that $(\tilde X,o)$ must therefore be Gorenstein. But,
there is no useful general result limiting the topology of complete
intersection singularities among the Gorenstein ones. So the statement
that a topologically defined Gorenstein singularity is a complete
intersection must be surprising (as also the very special form:
Brieskorn complete intersection).  There is one simple case for the
theorem: the universal abelian cover of a quotient singularity
$\C^2/G$ (with $G\subset \GL(2,\C)$ finite, acting freely off $(0,0)$)
is clearly $\C^2/G'$, where $G'$ is the commutator subgroup.  But then
$G'\subset \SL(2,\C)$, so $\C^2/G'$ must be smooth or a rational
double point, and necessarily a hypersurface.  In fact it is of type
$A_n$, $D_4$, $E_6$, or $E_8$, the only rational double points which
are Brieskorn singularities.

%From the topological point of view, the most natural class of $(X,o)$
%to consider for this conjecture are
The next impetus towards Conjecture \ref{conj:main} was our study of
normal surface singularities whose links are integral homology
spheres, developing out of \cite{neumann-wahl90}.  A $\Q$-Gorenstein
singularity of this type is necessarily Gorenstein, but the
conjecture would imply that it is always a complete intersection!
This issue was discussed by us in \cite{neumann-wahl90}, and will be
examined more closely in a forthcoming paper \cite{neumann-wahl00}.
There are a number of reasons to believe the conjecture in this case,
and this would imply completely unexpected topological restrictions on
links of such singularities.

Conjecture \ref{conj:main} applies to rational surface singularities,
but even the weaker result that every rational surface singularity is
the quotient of some complete intersection singularity is unknown.

Note finally that an equisingular deformation of a singularity
satisfying the Conjecture (e.g., a deformation of a weighted
homogeneous singularity as in Theorem \ref{th:weighted homog}) does
not itself obviously satisfy the conjecture! One problem is that if
$(Y,o)\to(X,o)$ is a finite covering (unramified off $o$), an
equisingular deformation of $X$ gives rise to a deformation of
$Y-\{o\}$, but not {\it a priori}\/ to one of $Y$.

\vspace{6pt}From the point of view of classification of
singularities, the easiest examples not covered by Theorem
\ref{th:weighted homog} are the \emph{quotient-cusps}. These are
singularities whose resolution graphs have the form
$$
\xymatrix@R=6pt@C=24pt@M=0pt@W=0pt@H=0pt{
  \overtag{\Circ}{-2}{8pt}\lineto[dr] && &&&
  \overtag{\Circ}{-2}{8pt}\lineto[dl]\\
  &\overtag{\Circ}{-e_1}{8pt}\lineto[r]
  &\overtag{\Circ}{-e_2}{8pt}\dashto[r]&\dashto[r]&
  \overtag{\Circ}{-e_k}{8pt}&&k\ge2,~~e_i\ge2, \text{ some }e_j>2.\\
  \overtag{\Circ}{-2}{8pt}\lineto[ur] && &&&
  \overtag{\Circ}{-2}{8pt}\lineto[ul]}$$
They are rational singularities, even log-canonical.  Like the
cusps, they are {\it taut}\/, i.e, the topology of the link
determines the analytic type of the singularity \cite {laufer}.
(See, e.g., \cite{kawamata} for a discussion of log-canonicals,
and \cite{hirzebruch}, \cite{laufer}, \cite{karras},
\cite{pinkham} and \cite{neumann81} for discussions of cusps and
quotient-cusps. The only log-canonical singularities which are not
weighted homogeneous are the cusps and quotient-cusps.)

The above quotient-cusp is double-covered by the cusp singularity
whose resolution graph is
$$
\xymatrix@R=8pt@C=24pt@M=0pt@W=0pt@H=0pt{
  &\overtag{\Circ}{-e_2}{8pt}\dashto[r]
  &&&\overtag{\Circ}{-e_{k-1}}{8pt}\dashto[l]\\
  \lefttag{\Circ}{2-2e_1}{8pt}\lineto[ur]\lineto[dr] &&&&&
  \righttag{\Circ}{2-2e_k}{8pt}\lineto[ul]\lineto[dl]\\
  &\overtag{\Circ}{-e_2}{8pt}\dashto[r]
  &&&\overtag{\Circ}{-e_{k-1}}{8pt}\dashto[l]}$$
It follows that the
universal abelian cover is also a cusp.  While it is easy to
determine, given the resolution graph, when a cusp is a complete
intersection (Proposition \ref{prop:complete} below), there is still
significant work to prove the surprising (but consistent): \emph{The
  universal abelian cover of a quotient-cusp is a complete
  intersection of embedding dimension $4$.}

To describe the precise result we need some notation.  We associate
to the quotient-cusp singularity described above the matrix
$$B=
\left(\begin{matrix}a&b\\c&d
  \end{matrix}
\right)%=M(e_1,\dots,e_k)
=B(e_1-1,e_2,\dots,e_{k-1},e_k-1)$$
with $$
B(b_1,\dots,b_k):= \left(
\begin{matrix}0&1\\-1&0
\end{matrix}\right)
\left(
\begin{matrix}0&-1\\1&b_k
\end{matrix}\right)
\dots
\left(
\begin{matrix}0&-1\\1&b_1
\end{matrix}\right)
.
$$
We shall see below (Lemma \ref{le:matrix}) that $B$ is an integer
unimodular matrix with positive entries, and each unimodular integer
matrix with positive entries arises this way for unique
$e_1,\dots,e_k$ as above.
Reversing the order of $e_1$ to $e_k$ just exchanges the entries $a$
and $d$.   Thus isomorphism types of quotient-cusp singularities are
in one-one correspondence with unimodular $2\times2$ matrices
{\scriptsize$\left(\begin{matrix}a&b\\ c&d
  \end{matrix}
\right)$} with
positive integer entries $a,b,c,d$, up to exchanging $a$ and $d$.

%The following is the precise version of our main Theorem.

\begin{theorem}\label{th:main} If $(V,p)$ is the quotient-cusp
  classified as above by {\scriptsize$\left(
    \begin{matrix}
      a&b\\ c&d
    \end{matrix}
\right)$}, then its
  universal abelian cover branched at $p$ is the complete intersection
  cusp singularity given by equations
$$xy=u^{2a}+v^{2a};\qquad uv=x^{2d}+y^{2d}$$
The covering degree is $16b$.
\end{theorem}

The cusp singularity in question is the one with resolution
graph
$$
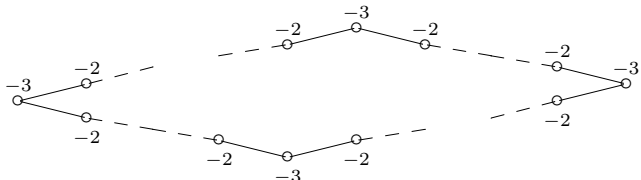
\xymatrix@R=2pt@C=24pt@M=0pt@W=0pt@H=0pt{
  &&&&&\overtag{\Circ}{-3}{6pt}\\
  &&&&\overtag{\Circ}{-2}{6pt}\lineto[ur] &&\overtag{\Circ}{-2}{6pt}\lineto[ul]\\
  &&&\dashto[ur] &&&&\dashto[ul]\\
  && &&&&&&\overtag{\Circ}{-2}{6pt}\dashto[ul]\\
  &\overtag{\Circ}{-2}{6pt}\dashto[ur]
  &&&&&&&&\overtag{\Circ}{-3}{6pt}
  \lineto[ul]\\
  \overtag{\Circ}{-3}{6pt}\lineto[ur]
  &&&&&&&&\undertag{\Circ}{-2}{3pt}\lineto[ur]\\
  &\undertag{\Circ}{-2}{3pt}\lineto[ul] &&&&&&\dashto[ur]\\
  &&\dashto[ul] &&&&\\
  &&&\undertag{\Circ}{-2}{3pt}\dashto[ul]
  &&\undertag{\Circ}{-2}{3pt}\dashto[ur]\\
  &&&&\undertag{\Circ}{-3}{3pt}\lineto[ul]\lineto[ur]\\ \\ \\}$$
where
the four strings of $-2$'s are lengths $2a-3$, $2d-3$, $2a-3$, and
$2d-3$. If $d=1$ or $a=1$ this must be replaced by
$$
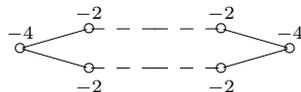
\xymatrix@R=3pt@C=24pt@M=0pt@W=0pt@H=0pt{
  &\overtag{\Circ}{-2}{6pt}\dashto[r]&\dashto[r]&\overtag{\Circ}{-2}{6pt}\\
  \overtag{\Circ}{-4}{6pt}\lineto[ur]\lineto[dr]&&&&
  \overtag{\Circ}{-4}{6pt}\lineto[ul]\lineto[dl]\\
  &\undertag{\Circ}{-2}{3pt}\dashto[r]&\dashto[r]&\undertag{\Circ}{-2}{3pt}\\
  \\ \\}$$
where the top and bottom strings are length $2a-3$ or
$2d-3$. The above equations for these cusps were given by Karras
\cite{karras}.

The point of our proof is that a cusp or quotient-cusp is determined
by its topology, so the proof becomes a purely topological (in fact
group-theoretical) calculation.

A quotient-cusp turns out to have other abelian covers which are
complete intersections, but we have no ``reason'' for it, except
calculation.  Also, any cusp singularity has a finite cover by a
complete intersection, but there are cusps with no abelian cover by a
complete intersection.

In the final section we sketch a more geometric proof of the theorem
by giving an explicit diagonal group action on a (necessarily
different) complete intersection representation of our cusp.  This is
in the spirit of the proof in \cite{neumann83b} of Theorem
\ref{th:weighted homog}, and borrows from themes we develop further in
\cite{neumann-wahl00}. But the verification that the quotient of the
action is correct is more difficult than in \cite{neumann83b}. So, as
presented here, this proof still depends on our earlier computations, and
for simplicity we only give full details for $b$ odd.

\section{Review of cusps and quotient-cusps}

In this section we review basic facts about cusp and quotient-cusp
singularities. Before describing the topological significance of the
classifying matrix for the quotient-cusp we prove the properties
mentioned above.  The above matrix
$$
B(b_1,\dots,b_k):= \left(
\begin{matrix}0&1\\-1&0
\end{matrix}\right)
\left(
\begin{matrix}0&-1\\1&b_k
\end{matrix}\right)
\dots
\left(
\begin{matrix}0&-1\\1&b_1
\end{matrix}\right)
.
$$
is clearly a unimodular integer matrix.

\begin{lemma}\label{le:matrix}
If $e_i\ge 2$ for each $i=1,\dots,k$ and at least one $e_i\ge3$,
then $$B(e_1-1,e_2,\dots,e_{k-1},e_k-1)$$  has positive entries.
Conversely, any unimodular integer matrix with positive entries
has a unique representation in the
above form.
\end{lemma}
\begin{proof}
Consider a product of matrices
$$\left(
  \begin{matrix}
    -\gamma&-\delta\\\alpha&\beta
  \end{matrix}
\right)=\left(
\begin{matrix}0&-1\\1&b_k
\end{matrix}\right)
\dots
\left(
\begin{matrix}0&-1\\1&b_1
\end{matrix}\right)
$$
with $k>1$.  The first statement of the lemma can be reformulated
that if $b_1\ge 1$, $b_i\ge2$ for $1<i\le k$ and $b_k\ge 1$, with at
least one of these inequalities strict, then
$\alpha,\beta,\gamma,\delta$ are all positive.

We will show, in fact, that $\alpha,\beta,\gamma,\delta$ are all
positive with the largest given by the following table:
\begin{align*}
\beta \text{~~largest}\quad\Longleftrightarrow\quad&b_1>1\text{ and } b_k>1\\
\alpha \text{~~largest}\quad\Longleftrightarrow\quad&b_1=1<b_k\\
\delta \text{~~largest}\quad\Longleftrightarrow\quad&b_1>1=b_k\\
\gamma \text{~~largest}\quad\Longleftrightarrow\quad&b_1=1=b_k
\end{align*}

The cases $(b_1,\dots,b_k)=(1,b)$, $(b,1)$, or $(1,b,1)$ can be
quickly verified by hand; we leave this to the reader.

Next suppose $b_1\ge2$ and $b_k\ge2$.  A simple induction on $k$ shows
that the entries of the product then satisfy $0< \gamma$,
$\gamma<\alpha<\beta$, and $\gamma<\delta<\beta$. This proves the
first case of the above table. Moreover, the same
induction shows that $\beta+\gamma-\alpha-\delta>0$ unless
$b_1=\dots=b_k=2$.
By post- and/or pre-multiplying by
{\scriptsize$ \left(
\begin{matrix}
0&-1\\1&1
\end{matrix}
\right)$} one now easily completes the proofs of the other three cases
of the table.

For the converse, consider
{\scriptsize$ \left(
\begin{matrix}
 -\gamma&-\delta\\ \alpha&\beta
\end{matrix}
\right)$}
with $\alpha,\beta,\gamma,\delta>0$ and $\alpha\delta-\beta\gamma=1$.
Assume first that $\beta$ is the largest of
$\alpha,\beta,\gamma,\delta$.  Then we must have $\gamma<\alpha<\beta$
and $\gamma<\delta<\beta$. Put
$b_1=\lceil\frac\delta\gamma\rceil$. If we write
$$\left(
\begin{matrix}-\gamma&-\delta\\\alpha&\beta
\end{matrix}\right)
=\left(
\begin{matrix}-\gamma'&-\delta'\\\alpha'&\beta'
\end{matrix}\right)\left(
\begin{matrix}0&-1\\1&b_1
\end{matrix}\right)
,
$$
then $\alpha'=b_1\alpha-\beta$, $\beta'=\alpha$,
$\gamma'=b_1\gamma-\delta$ and $\delta'=\gamma$. Thus, if $\gamma\ne
1$ we have $\gamma'>0$. Also,
$\alpha'>\frac\delta\gamma\alpha-\beta=\frac1\gamma>0$, and
$\delta'<\beta'$ and $\alpha'=b_1\alpha-\beta<(1+\frac\delta\gamma-
\frac1{\alpha\gamma})\alpha -\beta=\beta'$. It follows that
$\alpha',\beta',\gamma',\delta'$ are all positive and $\beta'$ is the
largest.  Thus we can proceed inductively. The value of $\gamma$ is
easily seen to decrease at each inductive step, and when $\gamma=1$ we
have $\gamma'=0$, so we can write
$$\left(
\begin{matrix}-\gamma'&-\delta'\\\alpha'&\beta'
\end{matrix}\right)=\left(
\begin{matrix}0&-1\\1&b_k
\end{matrix}\right).
$$  Thus, our matrix
{\scriptsize$\left(
  \begin{matrix}
    -\gamma&-\delta\\ \alpha&\beta
  \end{matrix}
\right)$}
is a product of the desired type with all $b_i\ge 2$.

The other three cases of the above table now follow easily as before
by first post- and/or pre-multiplying {\scriptsize$ \left(
\begin{matrix}
 -\gamma&-\delta\\ \alpha&\beta
\end{matrix}
\right)$}
by the inverse of {\scriptsize$ \left(
\begin{matrix}
0&-1\\1&1
\end{matrix}
\right)$}.

Finally, the uniqueness of the representation follows by noting that
the value we chose for $b_1$ in the above argument is forced by the
fact that $\gamma'\ge0$ by the calculation of the first part of the
proof.
\end{proof}

Let $M$ denote the oriented 3-manifold that is the total space of the
unit tangent bundle of the M\"obius band. This bundle has a section
consisting of unit tangent vectors that are parallel to the core
circle of the M\"obius band (once this circle is oriented); this also
provides a section on the boundary $\partial
M$. Thus we have a natural identification $\partial M\cong S^1\times
S^1$ where the first $S^1$ is base and the second $S^1$ is fiber.

\begin{lemma}\label{le:quotient-cusp}
  The link $N$ of the quotient-cusp with resolution graph
  $$
  
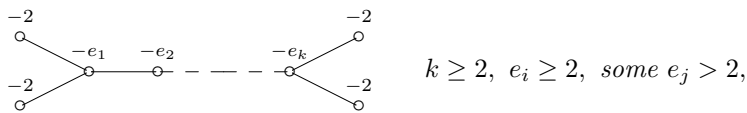
\xymatrix@R=6pt@C=24pt@M=0pt@W=0pt@H=0pt{
    \overtag{\Circ}{-2}{8pt}\lineto[dr] && &&&
    \overtag{\Circ}{-2}{8pt}\lineto[dl]\\
    &\overtag{\Circ}{-e_1}{8pt}\lineto[r]
    &\overtag{\Circ}{-e_2}{8pt}\dashto[r]&\dashto[r]&
    \overtag{\Circ}{-e_k}{8pt}&&k\ge2,~~e_i\ge2, \text{ some }e_j>2,\\
    \overtag{\Circ}{-2}{8pt}\lineto[ur] && &&&
    \overtag{\Circ}{-2}{8pt}\lineto[ul]}$$
  is the 3-manifold obtained
  as $N=M\cup_\phi M$, where $\phi$ is the pasting map $\phi\colon
  \partial M=T^2\to T^2=\partial M$ given by the matrix
  $B=B(e_1-1,e_2,\dots,e_{k-1},e_k-1)$ and $N$ is oriented with the
  orientation of the first $M$.

  Conversely, if the link of a singularity is diffeomorphic to
  $M\cup_\phi M$ for some $\phi\colon T^2\to T^2$ given by a unimodular
  matrix with positive entries, then the singularity is the
  corresponding quotient-cusp. In fact, it suffices that the link have
  the correct fundamental group.
\end{lemma}
  Since $\phi$ is orientation preserving in this lemma, $N$
  inherits opposite orientations from its two pieces $M$, which is why
  we must specify that we take the orientation inherited from the
  first.
\begin{proof}
As described in \cite{neumann81} (see ``Step 3'' on p.\ 314),
$N$ is obtained by plumbing from the plumbing diagram
$$
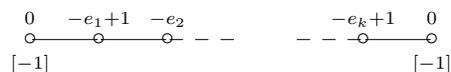
\xymatrix@R=12pt@C=24pt@M=0pt@W=0pt@H=0pt{\\
\overtag{\undertag{\Circ}{[-1]}{4pt}}{0}{8pt}\lineto[r]&
\overtag{\Circ}{-e_1+1}{8pt}\lineto[r]&
\overtag{\Circ}{-e_2}{8pt}\dashto[r]&&\dashto[r]&
\overtag{\Circ}{-e_k+1}{8pt}\lineto[r]&
\overtag{\undertag{\Circ}{[-1]}{4pt}}{0}{8pt}\\&}
$$
The two end vertices of this plumbing diagram give the pieces $M$ that
are pasted. As computed on page 319 of \cite{neumann81}, the
pasting map is given by the matrix product $JH_kJ\dots JH_1J$, where
$J=${\scriptsize$\left(
  \begin{matrix}
    0&1\\1&0
  \end{matrix}\right)$} and $H_i=${\scriptsize$
\left(
  \begin{matrix}
    -1&0\\b_i&1
  \end{matrix}\right)$} with $b_i=e_i$ or $e_i-1$ according as
$1<i<k$ or $i=1,k$. This gives the claimed gluing matrix after
left-multiplying by {\scriptsize${\left(
  \begin{matrix}
    1&0\\0&-1
  \end{matrix}\right)}$}
to make the gluing map unimodular (this reverses the orientation of
the second $M$).

The final part of the lemma says quotient-cusps are taut, in the sense
that a singularity with the same topology as a quotient-cusp \emph{is}
analytically that quotient-cusp. This was
proved by Laufer \cite{laufer}. The fact that the fundamental group
suffices to determine the topology is part of Corollary 4 of
\cite{neumann81}.
\end{proof}

A \emph{cusp singularity} is one whose minimal resolution graph is a
cycle of vertices representing rational curves $E_i$. We shall denote
this graph by $[-e_1,\dots,-e_k]$ (this notation is well-defined up to
cyclic permutation and reversal), where $-e_i$ is the Euler number of
the normal bundle of $E_i$. If $k>1$ then $-e_i$ is the
self-intersection number $E_i\cdot E_i$; but if $k=1$ then $E_1\cdot
E_1=-e_1+2$, since $E_1$ intersects itself normally in one point. The
$e_i$ satisfy $e_i\ge 2$ for all $i$ and some $e_j\ge 3$.
\begin{lemma}\label{le:cusp}
  The cusp singularity with resolution graph $[-e_1,\dots,-e_k]$ has
  link the torus bundle over $S^1$ with monodromy
$$A=\left(
\begin{matrix}0&-1\\1&e_k
\end{matrix}\right)
\dots
\left(
\begin{matrix}0&-1\\1&e_1
\end{matrix}\right).
$$
This matrix $A$ has trace $\ge3$.  The fundamental group of this
singularity link is the semi-direct product $\Z^2\rtimes\Z$ where a
generator of $\Z$ acts on $\Z^2$ by the matrix $A$.

Conversely, every torus bundle over $S^1$ with monodromy matrix $A$
having trace $\ge3$ is the link of a singularity: such a
matrix is conjugate to a matrix as above and the singularity is the
corresponding cusp.
\end{lemma}
\begin{proof} See \cite{neumann81}.
  The second part of the last sentence is again Laufer's tautness
  of cusp singularities \cite{laufer}. Our orientation convention here
  is different from \cite{neumann81} (where it was not made explicit),
  so the matrix given there was $JAJ$ rather than $A$.
\end{proof}

Although we will not need it, it is worth remarking that the cusp
double cover of the quotient-cusp with classifying matrix
{\scriptsize$\left(
  \begin{matrix}
    a&b\\ c&d
  \end{matrix}
\right)$} is the cusp with monodromy matrix
{\scriptsize$\left(
  \begin{matrix}
    d&b\\ c&a
  \end{matrix}
\right)\left(
  \begin{matrix}
    a&b\\ c&d
  \end{matrix}
\right)$}.

If one reverses the orientation of a torus bundle over the circle with
monodromy $A$ one gets a torus bundle with monodromy $A^{-1}$. Since
this does not change the trace of $A$, Lemma \ref{le:cusp} implies
that if one starts with the link of a cusp singularity then the result
is the link of another cusp singularity.  This is called the
\emph{dual cusp}. The following is well-known (e.g.,
\cite{hirzebruch}, \cite{neumann81}).
\begin{lemma}\label{le:dual cusp}
  To find the resolution graph of the dual cusp for a given cusp,
  replace each vertex weight $-e$ with $e\ge3$ in the resolution graph
  by a string of $e-3$ vertices weighted $-2$, and replace each
  intervening string of $d$ vertices weighted $-2$ ($d\ge0$) by a
  single vertex weighted $-d-3$.\qed
\end{lemma}

In particular, this lemma implies that the dual cusp to the cusp with
resolution cycle $[-e_1,\dots,-e_k]$  has resolution cycle of length
$$\sum_{i=1}^n(e_i-2).$$
For example, the resolution cycle for the
dual cusp to the cusp of Theorem \ref{th:main} is $[-2a,-2d,-2a,-2d]$
of length $4$.

This is useful in understanding which cusps are complete
intersections.  U. Karras proved:
\begin{proposition}[\cite{karras1}]\label{prop:complete}
  The cusp with resolution graph $[-e_1,\dots,-e_n]$ is a complete
  intersection if and only if $\sum_{i=1}^n(e_i-2)\le 4$. That is, the
  dual cusp has resolution cycle of length at most $4$.\qed
\end{proposition}

\section{Proof of Theorem \ref{th:main}}

Let $N$ and $M$ be as in Lemma \ref{le:quotient-cusp}.

 $\pi_1(M)=\langle x,y~|~x^{-1}yx=y^{-1}\rangle$, with $x$ representing a
core circle of the M\"obius band and $y$ representing the fiber. The
subgroup $\pi_1(\partial M)$ is generated by $x^2$ and $y$. The
van-Kampen theorem thus gives a presentation
\begin{align*}
\pi_1(N)=
\langle x_1,y_1,x_2,y_2~|&
~x_1^{-1}y_1x_1=y_1^{-1},\quad
x_2^{-1}y_2x_2=y_2^{-1},\\ &
\qquad\quad x_2^2=x_1^{2a}y_1^c,\quad
y_2 = x_1^{2b}y_1^d~~\rangle
\end{align*}
Abelianizing this presentation shows
\begin{align*}
\pi_1(N)^{ab}=
\langle X_1,Y_1,X_2,Y_2~|&
~2Y_1=0,\quad
2Y_2=0,\\ &
2X_2=2aX_1+cY_1,\quad
Y_2 = 2bX_1+dY_1~~\rangle
\end{align*}
The relation matrix
$$\left(
  \begin{matrix}
    0&2&0&0\\
0&0&0&2\\
2a&c&-2&0\\
2b&d&0&-1
  \end{matrix}\right)$$
of this abelian presentation has determinant $-16b$, so the order of
$\pi_1(N)^{ab}$ is $16b$.

The fundamental group of the universal abelian cover $\tilde N$ of $N$
is $$H:=\Ker\bigl(\pi_1(N)\to\pi_1(N)^{ab}\bigr).$$ We claim
\begin{lemma}
  $$ H=\langle s,x_1^{4b},y_1^2\rangle=\langle s,
x_2^{4b},y_2^2\rangle$$
where
$$s=(x_1x_2)^2x_1^{-2a-2}y_1^{-c}.$$
\end{lemma}
\begin{proof}
  Denote by $N_1$ the link of the cusp singularity that canonically
  double-covers our quotient-cusp. Since this is an abelian cover,
  $N_1$ is covered by the universal abelian cover $\tilde N$ of $N$.

  $N_1$ fibers over the circle $S^1$ with fiber the torus $T^2$,
  giving a semi-direct product representation
$$\pi_1(N_1)=V\rtimes C,$$
 where $V=\pi_1(T^2)$ is free abelian of rank $2$ and $C=\pi_1(S^1)$
 is infinite cyclic.

 $V$ is the abelian subgroup of $\pi_1(N)$ generated by $x_1^2$ and
 $y_1$. The elements $x_2^2$ and $y_2$ are related to these by the
 unimodular relations
$$
 \begin{matrix}
x_2^2&=&x_1^{2a}y_1^c\\ y_2 &= &x_1^{2b}y_1^d
 \end{matrix}
%x_2^2=x_1^{2a}y_1^c,\quad y_2 = x_1^{2b}y_1^d;
\qquad\text{with inverse}\qquad
 \begin{matrix}
%x_1^2=x_2^{2d}y_2^{-c},\quad y_1 = x_2^{-2b}y_2^a.
x_1^2&=&x_2^{2d}y_2^{-c}\\ y_1 &= &x_2^{-2b}y_2^a.
\end{matrix}
$$
 Thus $x_2^2$ and $y_2$ also generate $V$.  We can take $C$ to be
 generated by the element $$z=x_1x_2,$$ so the subgroup of $\pi_1(N)$
 generated by $z$ and $V$ is normal of index $2$.

 Now the subgroup $V_1$ of $V$ generated by $x_1^{4b}$ and $y_1^2$
 clearly has index $4b$ in $V$. Since $s$ equals $z^2$ modulo $V$, the
 subgroup generated by $s$ and $V_1$ has index $8b$ in $\pi_1(N_1)$,
 hence index $16b$ in $\pi_1(N)$.  On the other hand, it is evident
 from the presentation of $\pi_1(N)^{ab}$ that $s, x_1^{4b}$, and
 $y_1^2$ are in $H$. It follows that the subgroup they generate is $H$
 as claimed.

 Finally, the above relations show that $V_1$ is also generated by
 $x_2^{4b}$ and $y_2^2$, giving the second description of $H$ in the
 lemma.
\end{proof}

The above proof gives a semi-direct product representation
$$H=V_1\rtimes S,$$
where $S$ is infinite cyclic, generated by $s$.
The conjugation action of $s$ on $V_1$ is the same as the action of
$z^2$, so we first compute the action of $z$ on $V_1$. This will tell
us the covering $\hat N$ of $N_1$ determined by the subgroup $\langle
V_1,z\rangle$ of $\pi_1(N_1)$, and $\tilde N$ is then isomorphic to a
double covering of $\hat N$ ($\hat N$ is not necessarily itself a
covering of $N$). We have
\begin{gather*}
z^{-1}x_1^{4b}z=x_2^{-1}x_1^{-1}x_1^{4b}x_1x_2=
x_2^{-1}x_2^{4bd}y_2^{-2bc}x_2=
x_2^{4bd}y_2^{2bc}=\\
(x_1^{2a}y_1^{c})^{2bd}(x_1^{2b}y_1^d)^{2bc}=
(x_1^{4b})^{ad+bc}(y_1^2)^{2bcd}
\end{gather*}
\begin{gather*}
z^{-1}y_1^2z=x_2^{-1}y_1^{-2}x_2=x_2^{-1}x_2^{4b}y_2^{-2a}x_2
=x_2^{4b}y_2^{2a}=\\
(x_1^{2a}y_1^z)^{2b}(x_1^{2b}y_1^d)^{2a}=
(x_1^{4b})^{2a}(y_1^2)^{ad+bc}
\end{gather*}
Thus $z$ acts on $V_1$ by the matrix
$$\left(
  \begin{matrix}
    ad+bc&2a\\2bcd&ad+bc
  \end{matrix}\right)$$
This matrix is conjugate by
  {\scriptsize${\left(\begin{matrix}
    0&-1\\1&d
  \end{matrix}\right)}$} to the inverse of
$$
\left(
  \begin{matrix}
    0&-1\\1&2a
  \end{matrix}\right)
\left(
  \begin{matrix}
    0&-1\\1&2d
  \end{matrix}\right).
$$
So $\hat N$ is the orientation
reversal of the link of the cusp with resolution graph
$$\xymatrix@R=12pt@C=36pt@M=0pt@W=0pt@H=0pt{
  \lefttag{\Circ}{-2a}{6pt}\lineto@/^/[r]\lineto@/_/[r]&
  \righttag{\Circ}{-2d~.}{6pt}\\}$$
That is, $\hat N$ is the link
of the {dual cusp} to this cusp.  By Lemma \ref{le:dual cusp}
the dual cusp
in this case has graph ($a,d>1$)
$$
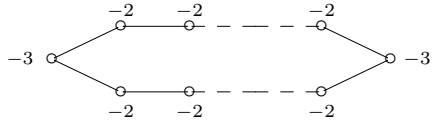
\xymatrix@R=8pt@C=24pt@M=0pt@W=0pt@H=0pt{
&\overtag{\Circ}{-2}{6pt}\lineto[r]&\overtag{\Circ}{-2}{6pt}\dashto[r]
&\dashto[r]&\overtag{\Circ}{-2}{6pt}\lineto[dr]\\
\lefttag{\Circ}{-3}{6pt}\lineto[ur]\lineto[dr]&&&&&
\righttag{\Circ}{-3}{6pt}\\
&\undertag{\Circ}{-2}{3pt}\lineto[r]&\undertag{\Circ}{-2}{3pt}\dashto[r]
&\dashto[r]&\undertag{\Circ}{-2}{3pt}\lineto[ur]\\ \\}$$
with the top and bottom strings of $-2$'s of length $2a-3$ and $2d-3$.
If $a$ or $d$ is $1$ the graph is
$$
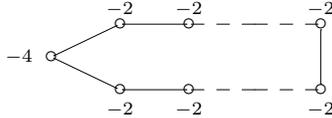
\xymatrix@R=8pt@C=24pt@M=0pt@W=0pt@H=0pt{
&\overtag{\Circ}{-2}{6pt}\lineto[r]&\overtag{\Circ}{-2}{6pt}\dashto[r]
&\dashto[r]&\overtag{\Circ}{-2}{6pt}\lineto[dd]\\
\lefttag{\Circ}{-4}{6pt}\lineto[ur]\lineto[dr]&&&&\\
&\undertag{\Circ}{-2}{3pt}\lineto[r]&\undertag{\Circ}{-2}{3pt}\dashto[r]
&\dashto[r]&\undertag{\Circ}{-2}{3pt}\\ \\}$$
with $2d-3$ or $2a-3$ vertices weighted $-2$.

Now the cusp we really want is the double cover whose monodromy acts
as the square of $z$ rather than $z$. It therefore has resolution
graph the double cover of the above cycle. This is as claimed after
Theorem \ref{th:main}, and the fact that this is a complete
intersection with the equations as described there can be found in
\cite{karras}.  This completes the proof.

Although we did not discuss orientations in detail, we were careful to
choose fundamental group generators appropriate to orientation in the
above proof.  This is important: the other orientation of our link
would mean we have the dual cusp, which is a complete intersection
only if $a+d=3$.  To be sure that orientations are correct we will
check a simple case by another method.

We first introduce terminology.  A link of a cusp singularity has a
fibration $N\to S^1$ with torus fibers. A covering of $N$ that pulls
back from a cyclic covering of the base $S^1$ will be called a
\emph{covering in the circle direction}. At the opposite extreme, a
\emph{fiberwise covering} will be one that induces the identity map on
the base circles of the torus fibrations.

Not every covering of the torus fiber of $N$ extends to a fiberwise
covering of $N$.  However the canonical $\Z/r\times\Z/r$ covering of
the fiber always does extend; if we write $\pi_1(N)$ as
$\Z^2\rtimes\Z$ then this cover corresponds to the subgroup
$r\Z^2\rtimes\Z$. Note that this subgroup depends on the choice of
splitting $\Z\to\pi_1(N)$ and is usually not a normal subgroup. Thus
the resulting fiberwise covering is not unique and not usually a
Galois covering. However, the resulting covering space is isomorphic
to the original link, since the action of $\Z$ on $r\Z^2$ is the same
as its action on $\Z^2$.

To confirm orientations in the above proof, consider the quotient-cusp
with matrix $A=${\scriptsize$\left(
  \begin{matrix}
    a&1\\ad-1&d
  \end{matrix}\right)$}.  The resolution graph in this case is
$$
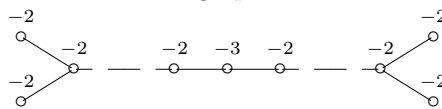
\xymatrix@R=8pt@C=18pt@M=0pt@W=0pt@H=0pt{
  \overtag{\Circ}{-2}{8pt}\lineto[dr] && &&&&&&
  \overtag{\Circ}{-2}{8pt}\lineto[dl]\\
  &\overtag{\Circ}{-2}{8pt}\dashto[r]&\dashto[r]
  &\overtag{\Circ}{-2}{8pt}\lineto[r]
  &\overtag{\Circ}{-3}{8pt}\lineto[r]
  &\overtag{\Circ}{-2}{8pt}\dashto[r]&\dashto[r] &
  \overtag{\Circ}{-2}{8pt}&&\\
  \overtag{\Circ}{-2}{8pt}\lineto[ur] && &&&&&&
  \overtag{\Circ}{-2}{8pt}\lineto[ul]}$$
with respectively $a-1$ and
$d-1$ vertices weighted $-2$ before and after the $-3$ on the central
string.  The universal abelian cover has degree 16 in this case and
the above proof shows it can be built up by the sequence of covers
corresponding to the decreasing chain of subgroups
$\pi_1(N)\supset\langle V,z\rangle\supset\langle V,s\rangle\supset
H=\langle V_1,s\rangle$. The first cover is the canonical two-fold
cover by a cusp. The resolution graph for this cusp has strings of
$2a-3$ and $2d-3$ weights $-2$ separated by $-3$'s.  Next is the
double cover of this cusp given by the subgroup $\langle
V,s\rangle\subset \langle V,z\rangle$. This is the double cover in the
circle direction, so it leads to a double cover of the resolution
graph.  The final cover is a $\Z/2\times\Z/2$ fiberwise cover of the
result, since $V/V_1=\Z/2\times\Z/2$.
This does not change the
resolution graph.  Thus, we have the right resolution graph.

\section{Other covers, covers of cusps}

In this section we make the following miscellaneous observations.

\begin{proposition}\label{prop:misc}
~\vspace{-1pt}
  \begin{enumerate}
  \item The quotient-cusp classified by
    {\scriptsize$\left(\begin{matrix} a&b\cr c&d
    \end{matrix}\right)$} has a natural $4b$-fold abelian cover
  giving the same complete intersection cusp as the $16b$-fold
  universal abelian cover.  The latter is a $\Z/2\times\Z/2$ fiberwise
  cover of the former.
  \item Every cusp can be covered by a complete intersection cusp ---
    even by a hypersurface cusp.
  \item Not every cusp has an abelian cover by a complete intersection
    cusp.
\end{enumerate}
\end{proposition}

It seems likely that there even exist cusps that have no Galois cover
by a complete intersection, but we have not checked this.

Before we prove this proposition we need to recall the ``discriminant
group'' of a singularity (for more details on the following, see e.g.,
\cite{looijenga-wahl}).

Given any isolated complex surface singularity $(V,o)$, let $X$ be a
tubular neighborhood of the exceptional divisor of a good resolution
so $N=\partial X$ is the link of the singularity. Then $H_2(X)\cong
\Z^n\subset H_2(X;\Q)\cong\Q^n$ where $n$ is the number of exceptional
curves in the resolution. Denote  by $S$ the
intersection form on these groups and define
$$H_n(X)^\#=\{v\in H_n(X;\Q):S(v,w)\in\Z\text{ for all }w\in
H_2(X)\}.$$
The embedding $H_2(X)\to H_2(X)^\#$ can be identified with
the map $H_2(X)\to H_2(X,N)$, so the long exact sequence in homology
identifies the \emph{discriminant group} $D:=H_2(X)^\#/H_2(X)$ with
the torsion of $H_1(N)$.  The intersection form $S$ induces on $D$ a
natural non-singular pairing
$$D\otimes D\to\Q/\Z, \quad v\otimes w\mapsto S(v,w)\text{ mod }\Z.$$
This is the torsion linking pairing of $N$.

If $K\subset D$
is a subgroup then we get an induced non-singular pairing
$$K\otimes (D/K^\perp)\to\Q/\Z$$
where $K^\perp$ is orthogonal
complement of $K$ under the pairing, so $D/K^\perp$ is canonically
isomorphic to the dual $\hat K=\Hom(K,\Q/\Z)$ and hence
non-canonically isomorphic to $K$ itself.

If $N$ is a rational homology sphere then the universal abelian cover
of $N$ is the Galois cover of $N$ determined by the natural
homomorphism $\pi_1(N)\to H_1(N)=D$. Thus any subgroup $K\subset D$
determines an abelian cover of $N$: the Galois cover with covering
transformation group $D/K$. The annihilator $K^\perp$ of $K$ gives a
cover that we call the \emph{dual cover} for $K$; its covering
transformation group is $D/K^\perp$.  The dual cover for $D$ is thus
the universal abelian cover.

For a cyclic subgroup generated by an element $[v]\in D$ the dual
cover can be described in terms of its extension to a branched cover
of the resolution $X$, branched along the exceptional divisor. If
$v=\sum\frac{c_i}{d_i}[E_i]$, where $[E_i]$ are the homology classes
of the exceptional curves $E_i$ of the resolution, then the branching
along $E_i$ is of degree $d/d_i$, where $d$ is the order of the $[v]$
(hence the degree of the covering).  The resulting branched cover of
$X$ may have quotient singularities at intersection points of the
$E_i$. After resolving these one obtains a (possibly non-minimal)
resolution of the cyclic covering of the given singularity. This will
be illustrated in the following.

\begin{proof}[Proof of Part 1 of Proposition \ref{prop:misc}]
  Using the notation just introduced, consider the following element
  $v\in D$ of order $2$ for a quotient-cusp of the Introduction
  (coefficients of $v$ are in parentheses and omitted if zero --- with
  vertices numbered $-1,0,1,\dots,k+2$ in the obvious way,
  $v=\frac12([E_{k+1}]+[E_{k+2}])$)
$$
\xymatrix@R=8pt@C=24pt@M=0pt@W=0pt@H=0pt{
  \overtag{\Circ}{-2}{8pt}\lineto[dr] && &&&
  \righttag{\overtag{\Circ}{-2}{8pt}}{(1/2)}{6pt}\lineto[dl]\\
  &\overtag{\Circ}{-e_1}{8pt}\lineto[r]
  &\overtag{\Circ}{-e_2}{8pt}\dashto[r]&\dashto[r]&
  \overtag{\Circ}{-e_k}{8pt}\\
  \overtag{\Circ}{-2}{8pt}\lineto[ur] && &&&
  \righttag{\overtag{\Circ}{-2}{8pt}}{(1/2)~.}{6pt}\lineto[ul]}$$
The double cover is
$$
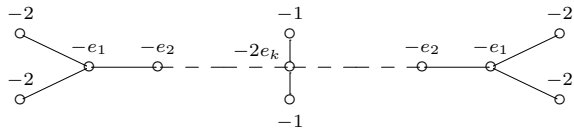
\xymatrix@R=8pt@C=24pt@M=0pt@W=0pt@H=0pt{
  \overtag{\Circ}{-2}{8pt}\lineto[dr] && &&\overtag{\Circ}{-1}{8pt}&
  &&&
  \overtag{\Circ}{-2}{8pt}\lineto[dl]\\
  &\overtag{\Circ}{-e_1}{8pt}\lineto[r]
  &\overtag{\Circ}{-e_2}{8pt}\dashto[r]&\dashto[r]^{-2e_k}&\Circ
%  \overtag{\Circ}{-2e_k}{8pt}
\dashto[r]\lineto[u]\lineto[d]&\dashto[r]
  &\overtag{\Circ}{-e_2}{8pt}\lineto[r]&\overtag{\Circ}{-e_1}{8pt}\\
  \overtag{\Circ}{-2}{8pt}\lineto[ur] && &&\undertag{\Circ}{-1}{4pt}&
  &&& \overtag{\Circ}{-2}{8pt}\lineto[ul]\\ \\}$$
which blows down to
$$
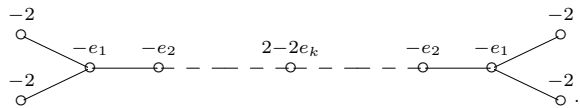
\xymatrix@R=8pt@C=24pt@M=0pt@W=0pt@H=0pt{
  \overtag{\Circ}{-2}{8pt}\lineto[dr] && &&& &&&
  \overtag{\Circ}{-2}{8pt}\lineto[dl]\\
  &\overtag{\Circ}{-e_1}{8pt}\lineto[r]
  &\overtag{\Circ}{-e_2}{8pt}\dashto[r]&\dashto[r]&
  \overtag{\Circ}{2-2e_k}{8pt} \dashto[r]&\dashto[r]
  &\overtag{\Circ}{-e_2}{8pt}\lineto[r]&\overtag{\Circ}{-e_1}{8pt}\\
  \overtag{\Circ}{-2}{8pt}\lineto[ur] && &&& &&&
  \righttag{\overtag{\Circ}{-2}{8pt}}{.}{6pt}\lineto[ul]\\}$$
Similarly, the element $w=\frac12([E_{-1}]+[E_{0}])$ gives a double
cover by a quotient-cusp, while the element $v+w$ gives the canonical
double cover by the cusp of the Introduction. These elements $v$ and
$w$ thus generate a Klein four-group $K\cong\Z/2\times\Z/2\subset D$
whose dual covering is the canonical covering by a cusp followed by
the canonical double covering that doubles the resolution cycle.
Denote this $4$-fold covering of the links $N_2\to N$.

Now this $K$ is self-orthogonal with respect to the linking form, so
$K\subset K^\perp$, so the dual covering $ N_3\to N$ determined by
$K^\perp$ factors through the covering just constructed. It is not
hard to check that $N_3\to N_2$ is a $\Z/b$-fold cover, but we do not
need this.  The universal abelian cover $\tilde N\to N$ factors
through $N_3$ and the covering transformation group for $\tilde N\to
N_3$ is $K$. Since this is $\Z/2\times \Z/2$ acting in the fibers of
$\tilde N$, the quotient is isomorphic to $\tilde N$ again. This
completes the proof of the first point of the above proposition.
\end{proof}

If the link $N$ of an isolated complex surface singularity $(V,o)$ is
not a rational homology sphere then there is no natural epimorphism
$\pi_1(N)\to D$ to the discriminant group, hence no natural Galois
cover with covering transformation group $D$.  However, if $(V,o)$ is
a cusp then the different epimorphisms of $H_1(N)=\Z\oplus D$ to $D$
are related by automorphisms of $\pi_1(N)$, and hence by automorphisms
of $(V,o)$ (see \cite{wall2}). So there is a natural cover up to
automorphisms, called the \emph{discriminant cover}. Hence for any
subgroup $K\subset D$ we may still speak of the \emph{cover for $K$}
and the \emph{dual cover for $K$}, meaning the covers with covering
transformation group $D/K$ respectively $D/K^\perp$.

The following proposition will be useful.
\begin{proposition}\label{prop:perp}
  If $D$ is the discriminant group for a cusp $(V,o)$ and $K\subset D$
  a subgroup, then the cover and the dual cover of $(V,o)$ for $K$ are
  mutually dual cusps.  In particular, taking $K=\{1\}$, then
  the discriminant cover of $(V,o)$ is the dual cusp of $(V,o)$.
\end{proposition}
\begin{proof}
  The link $N$ of the cusp is a torus bundle with monodromy $A$ say.
  We can consider $N$ to be $T^2\times [0,1]$ with the two boundary
  components identified by $A$. Replacing $T^2\times\{0\}\subset N$ by
  a tubular neighborhood and excising the interior of this
  neighborhood, we see $H_q(N,T^2)\cong H_{q}(T^2\times [0,1],
  T^2\times \{0,1\})\cong H_{q-1}(T^2)$. The long exact homology
  sequence for the pair $(N, T^2)$ thus leads to a long exact sequence
  $$\dots\to H_2(N)\to H_1(T^2)\buildrel{A-I}\over{\longrightarrow} H_1(T^2)\to
  H_1(N)\to \Ker(H_0(T^2)\to H_0(T^2))=\Z.$$
  This shows that the torsion $D$ of $H_1(N)$ is the
  image of $H_1(T^2)\to H_1(N)$ and is isomorphic to $\Z^2/(A-I)\Z^2$.

  If $C$ is a 1-cycle in $T^2$ then the 2-chain $C\times[0,1]\subset
  N$ has boundary $(A-I)C$.  It follows that the torsion
  linking form on $D=\Z^2/(A-I)\Z^2$ is given by
  $\ell([v],[w])=((A-I)^{-1}v).w$ (mod $\Z$), where $v,w\in\Z^2$ and
  $v.w$ is the standard skew-symmetric form on $\Z^2$ ($=$
  intersection form on $H_1(T^2)$).

  Suppose now that $K=W/(A-I)\Z^2\subset D$ and write
  $K^\perp=W'/(A-I)\Z^2$. Then lemma \ref{le:perp} below implies
  $W'=\frac1\Delta(A-I)^{-1}W$, where $\Delta=[D:K]/|\det(A-I)|$.

  The cover and the dual cover for $K$ of the given cusp have
  monodromy given by the action of $A$ on $W$ and $W'$ respectively.
  Let $w_1,w_2$ be an oriented basis of $W$ and put
  $w_1'=\frac1\Delta(A-I)^{-1}w_1,w_2'=\frac1\Delta(A-I)^{-1}w_2$. Then
  $w_1',w_2'$ is a reverse-oriented basis for $W'$, since $(A-I)^{-1}$
  has negative determinant (namely $\frac1{2-\tr(A)}$).  But the
  matrix of $A$ with respect to $w_1,w_2$ equals the matrix of $A$
  with respect to $w_1',w_2'$ since $A$ and $(A-I)^{-1}$ commute.
  Thus the two cusps are mutual duals.
\end{proof}

We used the following lemma, whose proof is an exercise.
  \begin{lemma}\label{le:perp}
Let $u.v$ denote the skew-symmetric form on
$\Q^2$ induced by the standard non-singular skew-symmetric form on $\Z^2$.
If $U\subset\Q^2$ is a $\Z$-lattice of covolume $\Delta$ (covolume is
  determinant of the matrix formed by an oriented basis) then the
  lattice
$$\{v:v.u\in\Z\text{ for all $u\in U\}$}$$
is $\frac1\Delta U$.\qed
  \end{lemma}

\begin{proof}[Proof of Parts 2,3 of Proposition \ref{prop:misc}]
Suppose we have a cusp given by the monodromy matrix
$A=${\scriptsize$\left(
  \begin{matrix}
    a&b\\c&d
  \end{matrix}
\right)$}. Assume $a\ne 0$ (if $a=0$ take $H=\Z^2$ in the following).
Then $A$ takes the subspace $H$ of $\Z^2$ generated by
{\scriptsize${\left(
  \begin{matrix}
    0\\1
  \end{matrix}
\right),\left(
  \begin{matrix}
    a\\c
  \end{matrix}
\right)}$}
to itself by the matrix
{\scriptsize${\left(
  \begin{matrix}
    0&-1\\1&t
  \end{matrix}
\right)}$}, where $t=\tr(A)=a+d$. Thus the cover defined by the
subgroup $H\rtimes\Z\subset\pi_1(N)$ is either the cusp with
resolution graph consisting of a cycle with one vertex weighted $-t$
or the dual cusp of this, according as the above basis is oriented
correctly or not, i.e., whether $a<0$ or $a>0$.  In any case, by going
to the discriminant cover if necessary (see Proposition
\ref{prop:perp}), we get the cusp with resolution graph dual to the
above length 1 cycle, so it is a cycle of $t-3$ $-2$'s and one $-3$.
This is a hypersurface cusp, proving the second statement of the
proposition. Since this is a fiberwise cover, this hypersurface cusp
not only covers, but is also covered by the original cusp. (In fact
any two cusps with the same trace of their monodromy are mutual
fiberwise covering spaces, see \cite{hirzebruch}.)

For part 3 of the Proposition, a necessary condition that a cusp
$(V,o)$ have no finite abelian cover by a complete intersection is
that
\begin{itemize}
\item neither $(V,o)$ nor the dual cusp $(V^*,o)$ is a complete
  intersection.
\end{itemize}
We shall see that if
\begin{itemize}
\item the discriminant group $D$ has prime order $p$ say,
\end{itemize}
then this necessary condition is also sufficient.  Indeed,
suppose that $(V,o)$ and $(V^*,o)$ are not complete intersections.  If
$D$ has prime order then any abelian cover is either a cover in the
circle direction or a cover in the circle direction of the
discriminant cover.  The resolution graph of the abelian cover is
therefore a cyclic covering of the resolution graph or the dual
resolution graph. By Proposition \ref{prop:complete} this cannot be the
resolution graph of a complete intersection.

An explicit example of this is the cusp with resolution cycle
$[2,4,2,2,5]$ with discriminant group of order $61$. This cusp is
self-dual and since the cycle is of length $5$ it is not a complete
intersection.
\end{proof}

\section{The explicit universal abelian cover}

Let $(V,p)$ be the quotient-cusp classified by {\scriptsize${\left(
  \begin{matrix}
    a&b\\c&d
  \end{matrix}
\right)}$}.
We will write down explicit equations and an explicit group action for
the universal abelian cover.

Let $\omega$ be a primitive $4b$-th root of $1$, and consider the
diagonal $4\times4$ matrices
\begin{align*}
  S_1&=[-\omega^a,\omega^a,\omega,\omega]\\
  S_2&=[\omega^a,-\omega^a,\omega,\omega]\\
  S_3&=[\omega,\omega,-\omega^d,\omega^d]\\
  S_4&=[\omega,\omega,\omega^d,-\omega^d].
\end{align*}
Let $G=\langle S_1,S_2,S_3,S_4\rangle$ be the abelian group they
generate.  We will show that $|G|=16b$, its group structure depending
on the parity of $c$.

\begin{theorem}
  Consider the complete intersection cusp singularity $(X,o)$ given by
  the equations
$$ x^2+y^2=u^\alpha v^\beta;\quad u^2+v^2=x^\gamma y^\delta,$$
where $\alpha, \beta, \gamma,\delta \ge0$ satisfy
$$\alpha+\beta=2a;\quad \gamma+\delta=2d;\quad
\alpha\equiv\beta\equiv\gamma\equiv\delta\equiv c \text{ \rm(mod $2$).}$$
Then
\begin{enumerate}
\item $G$ is an abelian group of order $16b$ which acts freely on a
punctured neighborhood of $o$ in $X$,
\item $(X/G, o)$ is isomorphic to $(V,p)$
\item $(X,o)\to(V,p)$ is the universal abelian cover of $(V,p)$.
\end{enumerate}
\end{theorem}
The proof shows that $G$ in fact acts freely on the smooth affine
variety $X-\{o\}$ unless $\alpha\beta=\gamma\delta=0$.
Note also that the equations for the cusp are \emph{not} in the
standard form of \cite{karras}; however, they are in the general form
we use in \cite{neumann-wahl00}, and should be considered as the
analog of ``Brieskorn complete intersections.''

We describe $G$ first.  Denote
$$ T=S_1S_2^{-1}=[-1,-1,1,1],$$
\begin{lemma}
  Every element of $G$ may be written uniquely as
$$S_1^jT^kS_3^l,\quad 0\le j \le 4b-1, 0\le k,l\le 1.$$
\end{lemma}
\begin{proof}
First note that $G$ contains $-I$, equal to $S_1^{2b}$ for $a$ odd, and
$S_1^{2b}T$ for $a$ even.  Since
$$S_2=TS_1,\quad S_4=-TS_3,\quad S_3^2=S_1^{2d}T^c,$$
every element of $G$ may be written as desired.  For uniqueness
consider
$$S_1^jT^kS_3^l=[(-1)^{j+k}\omega^{aj+l}, (-1)^k\omega^{aj+l},
(-1)^l\omega^{j+d}, \omega^{j+dl}],$$with $j,k,l$ as in the Lemma. The
quotient of the last two entries shows $l$ is determined.  The
quotient of the second term by the $a$-th power of the last then shows
$k$ is determined. This implies the uniqueness.
\end{proof}
\begin{lemma}\label{le:structure}
{\bf 1.}~~For $c$ even, $G\cong C_{4b}\times C_2\times C_2$, with factors
  generated by $S_1$, $T$, and $S_1^{-d}T^{c/2}S_3=[-1,1,-1,1]$.

The only non-identity elements of $G$ with $\ge2$ eigenvalues equal to
$1$ are the $6$ elements with $2$ $+1$'s and $2$ $-1$'s.

{\bf2.}~~For $c$ odd, $G\cong C_{4b}\times C_4$, with factors
generated by $S_1$ and $S_1^dS_3^{-1}=[(-1)^di^c,i^c,-1,1]$.

The only non-identity elements of $G$ with $\ge2$ eigenvalues equal to
$1$ are $\pm T$.
\end{lemma}
\begin{proof}
  The proof is straightforward and left to the reader. Note that in
  the second case $(S_1^dS_3^{-1})^2=T$.
\end{proof}
\begin{lemma}
  $G$ acts on $X$, freely off finitely many points (only $0$, unless
  $c$ is even and $\alpha\beta=\gamma\delta=0$).
\end{lemma}
\begin{proof}
  $G$ acts on $X$ because the generators $S_i$ each transform each
  defining equation by a character.  The origin is the only point of
  $X$ with at least $3$ coordinates equal to $0$.  A non-trivial element of $G$
  with $\ge2$ eigenvalues equal to $1$ must have fixed locus a
  coordinate plane given by setting two coordinates equal to $0$. Such
  a plane intersects $X$ in only finitely many points.
\end{proof}

The following facts are straightforward.  $X$ has one singular point,
at $0$.  The blowing-up $X_1$ of $X$ at $0$ is normal. The exceptional
fiber $E$ ($=$ projectivized tangent cone of $(X,o)$) consists of
either $4$ $\mathbb P^1$'s ($a,d>1$):
$$
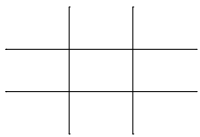
\xymatrix@R=16pt@C=24pt@M=0pt@W=0pt@H=0pt{&\lineto[ddd]&\lineto[ddd]\\
\lineto[rrr]&&&\\
\lineto[rrr]&&&\\
&&&}
$$
or $2$ $\mathbb P^1$'s ($a=1$ or $d=1$):
$$
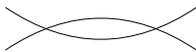
\xymatrix@R=16pt@C=24pt@M=0pt@W=0pt@H=0pt{
\lineto@/_12pt/[rrr]&&&\\
\lineto@/^12pt/[rrr]&&&\\&}
$$
In the first case, $X_1$ has singularities of type $A_{2a-3}$ and
$A_{2d-3}$ at alternate intersection points of the irreducible
exceptional curves. In the second case, $X_1$ has two $A_{2d-3}$
singularities. In particular, $(X,o)$ is the cusp singularity of
Theorem \ref{th:main}.

Since $G$ acts freely on the cusp $(X,o)$ off the origin, the
quotient $(X/G,o)$ must be either a cusp or a quotient-cusp.  But
the second (non-Gorenstein) case must occur; if $$\eta=dx\wedge
dy/(f_ug_v-f_vg_u)$$ is the standard nowhere-zero holomorphic
$2$-form on $X-\{o\}$  ($f$, $g$ the defining equations), one
finds
$$S_1^*\eta=-\eta.$$
So, $(X/G, o)$ is a quotient-cusp.
If we show $(X/G,o)$ has type {\scriptsize$\left(
  \begin{matrix}
    a&b\\c&d
  \end{matrix}
\right)$} then we are done, since the abelian cover $(X,o)\to (X/G,o)$
has degree $16b$, which is the degree of the universal abelian cover by
earlier computations.

We shall satisfy ourselves by carrying out this calculation for
$b=1$. This is not as special as it may seem; we show later that it
implies the general case that $b$ is odd.

So, let us assume $b=1$, and consider $a,d>1$ (the case that $a$ or
$d$ is $1$ is similar).  We will examine the quotient $X_1\to
X_1/G$. $G$ acts on the graph of $E$ as $C_2\times C_2$, either
rotating by $180^\circ$ or reflecting in a diagonal. So, the quotient
of $E$ by $G$ is one irreducible curve with $2$ distinguished
points. Further, off the $4$ singular points of $E$ the stabilizer of
any point is $\pm I$ (which acts trivially on $E$). This is because
$\pm I$ are the only scalar matrices on $G$. We conclude that on
$X_1/G$, the image of $E$ is a rational curve with two locally
irreducible
singular
points; and $X_1/G$ is smooth off these two points.
%$$
%\xymatrix@R=16pt@C=24pt@M=0pt@W=0pt@H=0pt{
%\lineto[rrr]&.&.&\\
%}$$
Further, the map $X_1\to X_1/G$ at a general point of $E$ consists
of a degree $2$ cover branched along $E$, followed by an $8$-to-$1$
covering map.

To find the singular points of $X_1/G$, choose a singular point of
$E$. Its stabilizer $H$ has order $8$. We do a local
calculation. Assume, e.g., $c$ is even. Then it is straightforward to
find local analytic coordinates so that $X_1$ is defined by $$
x_1^2+y_1^2=u^{2a-2},$$ $E$ is defined by $u=0$, and $H$ acts via
$[-i^{a-1}, i^{a-1}, i]$ and $[-1,-1,1]$. Writing down explicit
invariants, one computes that the quotient is a $D_{a+1}$-singularity;
further, the image of the curve $u=x_1\pm y_1=0$ will intersect the
resolved singularity on the long end. Therefore, the dual graph of the
minimal resolution $\tilde X$ of $X_1/G$ is
$$
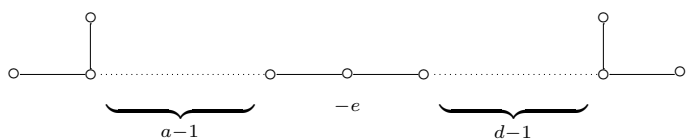
\xymatrix@R=16pt@C=24pt@M=0pt@W=0pt@H=0pt{
&\circ&&&&&&\circ\\
\circ\lineto[r]&\circ\lineto[u]\dotto[rr]&&\circ\lineto[r]
&\undertag\circ{-e}{3pt}\lineto[r]&\circ
\dotto[rr]&&\circ\lineto[r]\lineto[u]&\circ~,\\
&&\overtag{{}_{a-1}}{\hbox to 56pt{\upbracefill}}{6pt}&&&&
\overtag{{}_{d-1}}{\hbox to 56pt{\upbracefill}}{6pt}\\}
$$
where all but the central curves are $-2$-curves.

It remains only to show that $e=3$, as a calculation then shows that
the above is the  {\scriptsize${\left(
  \begin{matrix}
    a&1\\ad-1&d
  \end{matrix}
\right)}$} cusp quotient. This is achieved by computing the
divisor on $\tilde X$ of the function $g=x^4-u^4$ (this function
on $X$ is easily checked to be $G$-invariant).  On $X_1$, the
divisor of $g$ consists of $E$ counted with multiplicity $4$, plus
$16$ smooth curves, $4$ per component of $E$, intersection $E$
transversally at $4$ distinct non-singular points. By the
branching behavior of $X_1\to X_1/G$ along $E$, we conclude that
$g$ vanishes twice on $\tilde X$ along the central curve, and the
proper transform of $\{g=0\}$ consists of two smooth curves
intersecting the central curve transversally at distinct points. A
standard calculation then shows $e=3$.

We promised to show that the general case of $b$ odd follows from the
just completed case $b=1$ (the same argument will show that for
general $b$ it suffices prove the cases that $b$ is a power of
$2$).

We first note that it suffices
to show that $(X,o)\to(X/G, o)$ is a universal abelian cover.  For if
$(X/G,o)$ has type {\scriptsize${\left(
  \begin{matrix}
    a_1&b_1\\c_1&d_1
  \end{matrix}
\right)}$}, then the discriminant of $(X/G, o)$
would be
$$16b_1=|G|=16b$$
so $b=b_1$. Also, by Theorem \ref{th:main}, $\{a_1,d_1\}$  is
determined by the cusp cover, so, up to order, $a=a_1$, $d=d_1$.
Since $a_1d_1-b_1c_1=1$, we conclude that $(X/G, o)$ is of
type
\begin{scriptsize}$\left(
  \begin{matrix}
    a&b\\c&d
  \end{matrix}
\right)$\end{scriptsize}, as desired.

Now let $G_1\subset G$ be the subgroup of order $16$ generated by
$S=S_1^b$ and the other $1$ or $2$ generators described in Lemma
\ref{le:structure}. This is exactly the $G$ one would get if $b=1$
and, under the assumption that $b$ is odd, one easily checks that it
is the correct action.  We have already shown that $(X,o)\to
(X/G_1,o)$ is a universal abelian cover; this implies that the larger
abelian quotient $(X,o)\to(X/G,o)$ is also a universal abelian cover.

\end{document}